\documentclass[11pt]{article}

\usepackage{amsfonts}
\usepackage{amsmath,amssymb}

\def\RR{\mathbb{R}}
\def\NN{\mathbb{N}}
\def\EE{\mathbb{E}}
\def\AA{\mathcal{A}}

\def\II{{\mathbb{I}}}
\def\KK{\mathcal{C}}
\def\FF{\mathcal{F}}

\def\LL{\mathcal{L}}
\def\N{\mathcal{N}}

\def\a{\alpha}

\def\cl{{\rm cl}}

\def\proof{{\bf Proof. }}
\def\qed{\hfill$\square$\vspace{3mm}}

\def\t{\hat{\Theta}}

\newtheorem{theorem}{\bf Theorem}[section]
\newtheorem{lemma}[theorem]{\bf Lemma}
\newtheorem{definition}[theorem]{\bf Def{}inition}
\newtheorem{proposition}[theorem]{\bf Proposition}

\numberwithin{equation}{section}

\begin{document}

\renewcommand{\thefootnote}{}

\begin{center}
{\Large {\bf   
On symmetric fuzzy stochastic Volterra integral equations with retardation}}

\bigskip

Marek T. Malinowski

\bigskip

{\it Department of Applied Mathematics, \\ 
Tadeusz Ko\'sciuszko Cracow University of Technology,\\[0pt]
ul.\ Warszawska 24, 31-155 Krak\'ow, Poland\\
\bigskip
e-mail address: malinowskimarek@poczta.fm\\

\medskip
October 14, 2024}
\end{center}

\noindent
{\small {\bf Abstract:} This paper contains a study on stochastic Volterra integral equations with fuzzy sets-values and involving on a constant retardation. Moreover, the form of the equation is symmetric in the sense that fuzzy stochastic integrals are placed on both sides of the equation. We show that the considered initial value problem formulated in terms of symmetric fuzzy stochastic Volterra integral equation is well-posed. In particular, we show that there exists a unique solution and this solution depends continuously on the parameters of the equation. The results are achieved with the conditions of Lipschitz continuity of drift and diffusion coefficients, and continuity of kernels.}



\section{Preliminaries}

In this part of the paper we provide all the preliminary definitions, properties from multivalued analysis and fuzzy analysis, which constitute the background of the framework on fuzzy stochastic differential equations. Although all of them can be found in the author's earlier papers (cf.\ \cite{MAL12}), we also provide them here for the convenience of the reader.

By $\KK$ we denote
the set of all
nonempty, compact and convex subsets of $\RR^d$.
The distance between elements $A,B$ in $\KK$ can be measured by the Hausdorf{}f
metric $d_{H}$ described as 
$$
d_{H}\left( A,B\right) :=\max \left\{ \sup\limits_{a\in A}\inf\limits_{b\in
B}\|a-b\|,\sup\limits_{b\in B}\inf\limits_{a\in
A}\|a-b\|\right\},
$$
where $\|\cdot \|$ denotes a norm in $\RR^d$. It is known that the metric space $(\KK,d_H)$ is complete and
separable (see \cite{HU}). In $\KK$ one can define the addition and 
scalar multiplication in the following way for 
$A,B\in \KK$, $b\in\RR^d$, $\lambda\in\RR$ 
\begin{eqnarray}
&&A+B:=\{a+b: a\in A, b\in B\}, \,\,\, A+\{b\}:=\{a+b: a\in A\},\,\,\, \lambda A:=\{\lambda a: a\in A\}\nonumber.
\end{eqnarray}

A fuzzy set $u$ in $\RR^d$ is a generalization of an ordinary set (see \cite{ZAD,ZAD1}) and is characterized by its 
membership function $u\colon\RR^d\to[0,1]$. The value $u(x)$ (for each $x\in\RR^d$) 
is interpreted as the degree of membership of $x$ 
in the fuzzy set $u$. 

For a fuzzy set $u\colon\RR^d\to[0,1]$ by its $\a$-level, $\a\in(0,1]$, we mean the set $[u]^\a:=\{\,a\in\RR^d: u(a)\geqslant  \a \,\}$ and 
$[u]^0:=\cl\{\,a\in\RR^d: u(a)>0\,\}$ is called the support of $u$.

Let $\FF$ denote the family of fuzzy sets $u\colon\RR^d\to[0,1]$ such that
$[u]^\a\in\KK$ for every $\a\in[0,1]$.
Note that the set $\RR^d$ can be embedded into $\FF$ by the embedding $\langle\cdot\rangle\colon\RR^d\to\FF$ def{}ined as follows: for $r\in\RR^d$ we have 
$\langle r\rangle(x)=1$ if $x=r$, and $\langle r\rangle(x)=0$ if $x\not =r$.

Addition $u\oplus v$ and scalar multiplication $\beta\odot u$ in fuzzy set space $\FF$ can be def{}ined levelwise 
$$
[u\oplus v]^\a=[u]^\a+[v]^\a,\hspace{2mm}[\beta\odot u]^\a=\beta [u]^\a,
$$
where $u,v\in\FF$, $\beta\in\RR$ and $\a\in[0,1]$.

Let $u,v\in\FF$. If there exists $w\in\FF$ such that $u=v\oplus w$ then we call $w$ the Hukuhara dif{}ference of $u$ and $v$ and we denote it by $u\ominus v$. Note that $u\ominus v\not=u\oplus(-1)\odot v$. Also $u\ominus v$ may not exist, but if it exists it is unique. For $u,v\in\FF$ and $r_1,r_2\in\RR^d$ we have:
\begin{itemize}
\item[(P1)] $(u\oplus\langle r_1\rangle)\ominus\langle r_2\rangle=u\oplus\langle r_1 -r_2\rangle$,
\item[(P2)] the Hukuhara difference $(u\oplus \langle r_1\rangle)\ominus v$ exists if{}f $u\ominus v$ exists. Moreover, $(u\oplus \langle r_1\rangle)\ominus v=(u\ominus v)\oplus\langle r_1\rangle$.
\end{itemize}

A commonly used metric in $\FF$ is $d_\infty\colon\FF\times\FF\to[0,\infty)$ defined by
$$
d_\infty(u,v):=\sup\limits_{\alpha\in[0,1]}d_H([u]^\a,[v]^\a).
$$
It is known that $(\FF,d_\infty)$ is a complete metric space, but it is not separable and it is not locally compact. Useful properties of $d_\infty$ are listed below:\\
for every $u,v,w,z\in\FF$, $\beta\in\RR$ (see e.g.\ \cite{PUR})
\begin{itemize}
\item[(P3)] $d_\infty(u+w,v+w)=d_\infty(u,v)$,
\item[(P4)] $d_\infty(u+v,w+z)\leqslant  d_\infty(u,w)+d_\infty(v,z)$,
\item[(P5)] $d_\infty(\beta u, \beta v)=|\beta|d_\infty(u,v)$,
\item[(P6)] $d_\infty(u\ominus v,\langle0\rangle)=d_\infty(u,v)$,
\item[(P7)] $d_\infty(u\ominus v,u\ominus w)=d_\infty(v,w)$,
\item[(P8)] $d_\infty(u\ominus v,w\ominus z)\leqslant  d_\infty(u,w)+d_\infty(v,z)$.
\end{itemize}

Let
$\left( \Omega,\AA, P \right)$ be a complete probability space and
${\cal M}(\Omega,\AA;\KK)$
denote the family
of $\AA$-measurable multivalued mappings $F\colon\Omega\to\KK$ (multivalued random variable) such that
$$
\left\{ \omega\in \Omega : F(\omega)\cap O\neq \emptyset \right\}
\in \AA\hspace{2mm}\mbox{ for
every open set}\hspace{2mm} O\subset \RR^d.
$$
A multivalued random variable $F\in {\cal M}(\Omega,\AA; \KK)$
is called $L^{p}$-integrally 
bounded, $p\geqslant  1$, if there exists
$h\in L^{p}\left(\Omega,\AA, P;\RR\right)$
such that $\|a\|\leqslant  h(\omega)$ for any $a$ and $\omega$ with $a\in F(\omega)$. 
It is known (see \cite{HIA}) that $F$ is $L^p$-integrally  bounded if{}f
$\omega\mapsto d_{H}\left(F(\omega),\{0\}\right)$ is in $L^p(\Omega,\AA,P;\RR)$, 
where $L^p(\Omega,\AA, P;\RR)$ is a space of equivalence
classes (with respect to the equality $ P$-a.e.) of
$\AA$-measurable random variables $h\colon\Omega\to \RR$ such that
$\EE |h|^p=\int_\Omega |h|^p d P<\infty$. 
Let us denote
$$
\LL^{p}(\Omega,\AA,P;\KK):=
\left\{ F\in {\cal M}(\Omega,\AA;\KK) :
F\hspace{2mm}\mbox{is}\hspace{2mm}L^{p}\mbox{-integrally  bounded} \right\},\,\,\, p\geqslant 1.
$$
The multivalued random variables $F,G\in\LL^{p}\bigl(\Omega,\AA, P;\KK\bigr)$
are considered to be identical, if $F=G$ holds $P$-a.e.

Let $T>0$, and denote $I:= [0,T]$. Let the system
$(\Omega,\AA,\{\AA_{t}\}_{_{t\in I}}, P)$ be a
complete, f{}iltered probability space with a f{}iltration
$\{\AA_{t}\}_{_{t\in I}}$ satisfying the usual hypotheses,
i.e., $\{\AA_{t}\}_{_{t\in I}}$ is an increasing and right
continuous family of sub-$\sigma$-algebras of $\AA$, and $\AA_0$
contains all $P$-null sets.
We call $X\colon I\times\Omega\to \KK$ a multivalued
stochastic process, if for every $t\in  I$ a mapping
$X(t)\colon\Omega\to\KK$ is a multivalued random
variable.
We say that a multivalued stochastic process $X$ is $d_H$-continuous, if 
almost all (with respect to the probability measure $P$) 
its paths, i.e., the mappings 
$X(\cdot,\omega)\colon I\to\KK$ are the $d_H$-continuous functions.
A multivalued stochastic process $X$ is said to be 
$\{\AA_t\}_{t\in I}$-adapted, if for every $t\in I$ the multivalued random variable
$X(t)\colon\Omega\to \KK$ is $\AA_t$-measurable.
It is called measurable, if
$X\colon  I\times\Omega\to \KK$ is a
${\cal B}(I)\otimes\AA$-measurable multivalued random variable, 
where ${\cal B}(I)$ denotes the Borel
$\sigma$-algebra of subsets of $I$.
If $X\colon I\times\Omega\to \KK$  is
$\{\AA_t\}_{t\in I}$-adapted and measurable, then it will be
called nonanticipating.
Equivalently, $X$ is nonanticipating if{}f 
$X$ is measurable with respect to the $\sigma$-algebra $\N$
which is def{}ined as follows
$$
\N:=\{A\in {\cal B}(I)\otimes\AA :A^{t}\in\AA_{t}
\hspace{2mm}\mbox{for every}\hspace{2mm}
t\in  I\},
$$
where $A^{t}=\{\omega :(t,\omega )\in A\}$.
A multivalued nonanticipating stochastic process
$X\colon I\times\Omega\to\KK$ is called 
$\LL^p$-integrally bounded, if there exists a measurable 
stochastic process $h\colon I\times\Omega\to\RR$ such that $\EE\Bigl(\int_I|h(s)|^pds\Bigr)<\infty$ and   
$|||X(t,\omega)|||\leqslant  h(t,\omega)$ for a.a.\ $(t,\omega)\in  I\times\Omega$.
By $\LL^p(I\times\Omega,\N;\KK)$ we denote the set of all 
equivalence classes (with respect to the equality $\gamma\times P$-a.e., $\gamma$ denotes the Lebesgue measure) 
of nonanticipating and $\LL^p$-integrally bounded multivalued stochastic processes.

A mapping $x\colon\Omega\to\FF$ is said to be
a fuzzy random variable (see \cite{PUR}), if
$[x]^\a\colon\Omega\to \KK$ is an $\AA$-measurable multivalued random variable
for all $\a\in[0,1]$. It is known from \cite{KIM}
that $x\colon\Omega\to\FF(\RR^d)$ is the fuzzy random variable if{}f
$x\colon(\Omega,\AA)\to(\FF(\RR^d),{\cal B}_{d_S})$ is $\AA|{\cal B}_{d_S}$-measurable,
where $d_S$ denotes the Skorohod metric in $\FF(\RR^d)$ and ${\cal B}_{d_S}$ denotes the $\sigma$-algebra
generated by the topology induced by $d_S$. A fuzzy random variable $x\colon\Omega\to\FF$
is said to be $L^p$-integrally  bounded, $p\geqslant 1$,
if $[x]^0$ belongs to $\LL^p(\Omega,\AA, P;\KK)$.
By $\LL^p(\Omega,\AA, P;\FF)$ we denote the set of all 
$L^p$-integrally  bounded fuzzy random variables, where we consider
$x,y\in\LL^p(\Omega,\AA, P; \FF)$ as identical if
$x=y$ holds $P$-a.e.

We call $x\colon I\times\Omega\to \FF$ a fuzzy
stochastic process, if for every $t\in  I$ the mapping
$x(t,\cdot)\colon\Omega\to\FF$ is a fuzzy random
variable. 
We say that a fuzzy stochastic process $x$ is $d_\infty$-continuous, if 
almost all (with respect to the probability measure $P$) 
its trajectories, i.e.\ the mappings 
$x(\cdot,\omega)\colon I\to\FF$ are the $d_\infty$-continuous functions.
A fuzzy stochastic process $x$ is called 
$\{\AA_t\}_{t\in I}$-adapted, if for
every $\a\in[0,1]$ the multifunction
$[x(t)]^\a\colon\Omega\to \KK$ is $\AA_t$-measurable
for all $t\in I$.
It is called measurable, if
$[x]^\a\colon  I\times\Omega\to \KK$ is a
${\cal B}( I)\otimes\AA$-measurable multifunction
for all $\a\in[0,1]$, where ${\cal B}( I)$ denotes the Borel
$\sigma$-algebra of subsets of $ I$.
If $x\colon  I\times\Omega\to \FF$  is
$\{\AA_t\}_{t\in I}$-adapted and measurable, then it is
called nonanticipating.
Equivalently, $x$ is nonanticipating if{}f for every $\a\in[0,1]$ the
multivalued random variable $[x]^\a$ is measurable with respect to the $\sigma$-algebra $\N$. A fuzzy stochastic process
$x$ is called $L^p$-integrally  bounded ($p\geqslant  1$), if there exists a measurable stochastic process $h\colon I\times\Omega\to\RR$ such
that $\EE\int_{I}|h(t)|^pdt<\infty$ and $d_\infty(x(t,\omega),\t)\leqslant  h(t,\omega)$ for a.a.\ $(t,\omega)\in  I\times\Omega$.
By $\LL^p( I\times\Omega,\N;\FF)$ we denote
the set of nonanticipating and $L^p$-integrally  bounded fuzzy stochastic
processes.


Let $x\in\LL^p(  I\times\Omega,\N;\FF)$, $p\geqslant 1$.
For such the process $x$  
we can def{}ine (see e.g.\ \cite{MAL1,MAL2}) the fuzzy stochastic 
Lebesgue--Aumann integral which is a fuzzy random variable 
$$
\Omega\ni\omega\mapsto\int_Ix(s,\omega)ds\in\FF.
$$
Then $\int_0^tx(s)ds$ (from now on we do not write the argument $\omega$) 
is understood as $\int_I{\bf 1}_{[0,t]}(s)x(s)ds$.
For the fuzzy stochastic Lebesgue--Aumann integral we have the following properties (see \cite{MAL1,MAL2}).

\begin{proposition}\label{PRO}
Let $p\geqslant  1$. If $x,y\in\LL^p([0,T]\times\Omega,\N;\FF)$ 
then \begin{itemize}
\item[(i)] $[0,T]\times\Omega\ni(t,\omega)\mapsto\int_0^tx(s,\omega)ds\in\FF$ belongs to $\LL^p( [0,T]\times\Omega,\N;\FF)$,
\item[(ii)] the fuzzy stochastic  process $(t,\omega)\mapsto\int_0^tx(s,\omega)ds$ is $d_\infty$-continuous,
\item[(iii)] with probability one, for every $t\in[0,T]$
$$d_\infty^p\bigl(\int_0^t x(s)ds,\int_0^ty(s)ds\bigr)\leqslant t^{p-1}\int_0^t d_\infty^p(x(s),y(s))ds,$$
\item[(iv)] for every $t\in [0,T]$ it holds 
$$\EE \sup _{z\in [0,t]} d_\infty^p\bigl(\int_0^zx(s)ds,\int_0^zy(s)ds\bigr)\leqslant  t^{p-1}
\EE \int_0^td_\infty^p\bigl(x(s),y(s)\bigr)ds.
$$
\end{itemize}
\end{proposition}

As we mentioned e.g.\ in \cite{MAL1,MAL2} it is not possible to define fuzzy stochastic integral of It\^o type such the integral in such a fashion that it is not a crisp random variable. Hence, we will consider the diffusion part of the fuzzy stochastic differential equation as the crisp stochastic It\^o integral whose values are  embedded into $\FF$.

\section{Results concerning qualitative properties of solutions}

The main object of investigations presented in this paper is the symmetric fuzzy stochastic Volterra integral equation with retardation and this equation is of the following form
\begin{eqnarray}\label{glowne}
&&\hspace{-1.4cm}U(t)\oplus\int_0^tg_1(t,s) K_1(U(s),U(s-\tau))ds\oplus\bigg\langle\int_0^t h_1(t,s)L_1(U(s),U(s-\tau))dB_1(s)\bigg\rangle\nonumber\\
&&\hspace{-1.4cm}=\ \Phi(t)\oplus\int_0^tg_2(t,s) K_2(U(s),U(s-\tau))ds\oplus\bigg\langle\int_0^t h_2(t,s)L_2(U(s),U(s-\tau))dB_2(s)\bigg\rangle\\
&&\hspace{-1.4cm} \text{for } t\in[0,T],\quad\text{ and } U(t)=\Phi(t) \text{ for }t\in[-\tau,0],\nonumber
\end{eqnarray}
where $\tau>0$ symbolizes a constant retardation, $\Phi:[-\tau,T]\to\FF$ is a fixed fuzzy set-valued mapping, 
$K_1,K_2:\FF\times\FF\to\FF$ stand for the drift coefficients of the equation, and
$G_1,G_2:\FF\times\FF\to\RR$ denote the diffusion coefficients of the equation. The considered equation is driven by 
the real-valued $\{\mathcal{A}_t\}_{t\in[0,T]}$-Brownian motions (not necessarily independent) $B_1$, $B_2$, and contains the continuous kernels $g_1, g_2,h_1,h_2:[0,T]\times[0,T]\to\RR$. In this section, our aim is to present important results of the analysis of the qualitative properties of the solutions of Equation~\eqref{glowne}, i.e.\ we will justify the existence of the solution and its uniqueness under certain conditions imposed on the coefficients of the equation, and we will show that the solution does not change significantly when the coefficients of the equation change slightly.

Note that the proposed symmetric form of Equation~\eqref{glowne} cannot be equivalent to the one-sided standard form, i.e.\ only with integrals on the right-hand side. In fact, what is only possible is to rewrite this equation in the form
\begin{eqnarray*}
U(t)&=&\Bigl[\Phi(t)\oplus\int_0^tg_2(t,s) K_2(U(s),U(s-\tau))ds\Bigr]\ominus \int_0^tg_1(t,s) K_1(U(s),U(s-\tau))ds\\
&&\oplus\bigg\langle\int_0^t h_2(t,s)L_2(U(s),U(s-\tau))dB_2(s)-\int_0^t h_1(t,s)L_1(U(s),U(s-\tau))dB_1(s)\bigg\rangle.
\end{eqnarray*}
This form of Equation~\eqref{glowne} shows that there will be requirements for necessary assumptions about the existence of certain Hukuhara differences, and this will be done later.

We begin our analysis with a formal description of what the solution to Equation~\eqref{glowne} is.
Let $\tilde{T}\in (0,T]$. 
\begin{definition}\label{ROZ}
A solution to Equation~\eqref{glowne} on the interval $[-\tau,\tilde{T}]$ is a fuzzy stochastic process $U\colon [-\tau,\tilde{T}]\times\Omega\to \FF$ which satisfies:
\item (1) $U(t)=\Phi(t)$ for every $t\in[-\tau,0]$ with probability one,
\item (2) $U\big|_{[0,\tilde{T}]\times\Omega}\in\LL^2([0,\tilde{T}]\times\Omega,\N;\FF)$,
\item (3) $U$ is $d_\infty$-continuous, i.e. $P$-a.a.\ trajectories of $U$ are the $d_\infty$-continuous mappings,
\item (4) $U$ makes Equation~\eqref{glowne} true.
\end{definition}
From a practical point of view, the most favorable situation is when the equation has only one solution. It is known that stochastic processes can be identified in several different ways. Therefore, we also provide how we understand the uniqueness of the solution.
\begin{definition}\label{JED}
We say that a solution $U\colon [-\tau,\tilde{T}]\times\Omega\to  \FF$ to
 Equation~\eqref{glowne} is unique, if for any other solution $V\colon [-\tau,\tilde{T}]\times\Omega\to  \FF$ it holds
$$P(U(t)=V(t)\text{ for every }t\in[0,\tilde{T}])=1.$$
\end{definition}

Even classical single-valued stochastic differential equations may have no solutions, so the same may happen in the study of fuzzy stochastic differential or integral equations. Therefore, an important part of the research is proving the existence of a solution in connection with the conditions that the coefficients and parameters of the equation must meet. Below we formulate the requirements that we will impose to obtain the existence of a unique solution to Equation~\eqref{glowne}.
Apart from continuity of the kernels $g_1,h_1,g_2,h_2$ in the integrals, we will require that
\begin{itemize}
\item[(H{}0)] for the initial mapping $\Phi$ it holds $\sup\limits_{t\in[-\tau,T]}\|\Phi(t)\|_\FF<\infty$,
\item[(H{}1)]  for every $U_1,V_1,U_2,V_2\in\FF$
\begin{eqnarray*}
&&\hspace{-1cm}\max\big\{d_\infty^2(K_1(U_1,V_1),K_1(U_2,V_2)),
d_\infty^2(K_2(U_1,V_1),K_2(U_2,V_2)),|L_1(U_1,V_1)-L_1(U_2,V_2)|^2,\\
&& \hspace{1cm}|L_2(U_1,V_1)-L_2(U_2,V_2)|^2\big\}\leqslant C(d^2_\infty (U_1,U_2)+d^2_\infty (V_1,V_2)),
\end{eqnarray*}
where $C$ is a certain positive constant,
\item[(H{}2)] for every $U,V\in \FF$
\begin{eqnarray*}
\max\big\{\|F_1(U,V)\|^2_\FF,\|F_2(U,V)\|^2_\FF,|G_1(U,V)|^2,|G_2(U,V)|^2\big\}
\leqslant C(1+\|U\|^2_\FF+\|V\|^2_\FF),
\end{eqnarray*}
where $C$ is a certain positive constant,
\item[(H{}3)] there exists $\tilde{T}\in(0,T]$ such that for every $n\in\NN\cup\{0\}$ the mappings $U^n\colon \tilde{I}\times\Omega\to\FF$, where $\tilde{I}=[0,\tilde{T}]$, described as
$$
U^0(t)=\Phi(t)\text{ for }t\in[-\tau,\tilde{T}] \text{ with probability one, and for }n\in\NN
$$
\begin{align}\label{ninfty-picard-3}
    U^n(t)&=\Bigl[\Phi(t)\oplus\int_{0}^{t}g_2(t,s)K_{2}\left(U^{n-1}(s),U^{n-1}(s-\tau)\right)ds\Bigr]\\\notag
&\quad\ominus \int_{0}^{t}g_1(t,s)K_{1}\left(U^{n-1}(s),U^{n-1}(s-\tau)\right)ds\\\notag
       &\quad\oplus\Bigl\langle\int_{0}^{t}h_2(t,s)L_{2}(U^{n-1}(s),U^{n-1}(s-\tau))dB_2(s)\\\notag
&\quad-\int_{0}^{t}h_1(t,s)L_{1}\left(U^{n-1}(s),U^{n-1}(s-\tau)\right)dB_1(s)\Bigr\rangle\notag
\end{align}
are well defined (in particular the Hukuhara differences do exist). 
\end{itemize}
Using this set of requirements, we will justify the existence of a unique solution for symmetric fuzzy stochastic Volterra integral equation with retardation. First, however, we will start with a useful lemma dealing with the fuzzy stochastic processes described in (H3).

\begin{lemma}\label{OGR}
Suppose that (H0), (H2) and (H3) are fulfilled. Then, for the sequence $\{U^n\}_{n=0}^\infty$ defined in (H3) we have the property:
 $$\text{there is a positive constant }M\text{ such that for every }n\in\NN \text{ it holds}$$
$$\sup_{t\in[-\tau,\tilde{T}]}\EE \|U^n(t)\|^2_\FF\leqslant M.$$
\end{lemma}
\proof
Let us begin with an observation that for every $n\in\NN$ and every $t\in[-\tau,0]$ it holds $U^n(t)=\Phi(t)$. Therefore,  $\sup_{t\in[-\tau,0]}\EE \|U^n(t)\|^2_\FF=\sup_{t\in[-\tau,0]} \|\Phi(t)\|^2_\FF$ and the latter expression is bounded by hypothesis (H0). Continuing, let us consider $\sup_{z\in[0,t]}\EE \|U^{n}(t)\|^2_\FF$ with fixed $n\in\NN$ and fixed $t\in[0,\tilde{T}]$. Firstly, notice that
\begin{eqnarray*}
\sup_{z\in[0,t]}\EE \|U^{n}(t)\|^2_\FF
  &= &\sup_{z\in[0,t]}\EE \bigg\|\bigg[\Phi(z)\oplus\int_{0}^{z}g_2(z,s)K_{2}(U^{n-1}(s),U^{n-1}(s-\tau))ds\bigg]\\
&&\ominus \int_{0}^{z}g_1(z,s)K_{1}\left(U^{n-1}(s),U^{n-1}(s-\tau)\right)ds\\
  &&\oplus\bigg\langle\int_{0}^{z}h_2(z,s)L_{2}(U^{n-1}(s),U^{n-1}(s-\tau))dB_2(s)\\
&&-\int_{0}^{z}h_1(z,s)L_{1}(U^{n-1}(s),U^{n-1}(s-\tau))dB_1(s)\bigg\rangle\bigg\|^2_\FF.
 \end{eqnarray*}
Hence, applying the Cauchy--Schwarz inequality we arrive at
\begin{eqnarray*}
&&\sup_{z\in[0,t]}\EE \|U^{n}(t)\|^2_\FF\\
 &\leqslant&5\sup_{z\in[0,t]}\|\Phi(t)\|^2_\FF+5\sup_{z\in[0,t]}\EE \bigg\|\int_{0}^{z}g_2(z,s)K_{2}(U^{n-1}(s),U^{n-1}(s-\tau))ds\bigg\|^2_\FF\\
&&+\ 5\sup_{z\in[0,t]}\EE \bigg\|\int_{0}^{z}g_1(z,s)K_{1}(U^{n-1}(s),U^{n-1}(s-\tau))ds\bigg\|^2_\FF\\
&&+\ 5\sup_{z\in[0,t]}\EE \bigg|\int_{0}^{z}h_2(z,s)L_{2}(U^{n-1}(s),U^{n-1}(s-\tau))dB_2(s)\bigg|^2\\
&&+\ 5\sup_{z\in[0,t]}\EE \bigg|\int_{0}^{z}h_1(z,s)L_{1}(U^{n-1}(s),U^{n-1}(s-\tau))dB_1(s)\bigg|^2.
\end{eqnarray*}
Referring to Proposition 2.1 (iii), the It\^o isometry, continuity of the kernels and the Fubini Theorem
\begin{eqnarray*}
&&\sup_{z\in[0,t]}\EE \|U^{n}(t)\|^2_\FF\\
 &\leqslant&5\sup_{z\in[0,t]}\|\Phi(t)\|^2_\FF+5\sup_{z\in[0,t]} z\EE \int_{0}^{z}\|g_2(z,s)K_{2}(U^{n-1}(s),U^{n-1}(s-\tau))\|^2_\FF ds\\
&&+\ 5\sup_{z\in[0,t]}z\EE\int_0^z \|g_1(z,s)K_{1}(U^{n-1}(s),U^{n-1}(s-\tau))\|^2_\FF ds\\
&&+\ 5\sup_{z\in[0,t]}\EE\int_0^z \big|h_2(z,s)L_{2}(U^{n-1}(s),U^{n-1}(s-\tau))\big|^2ds\\
&&+\ 5\sup_{z\in[0,t]}\EE \int_0^z \big|h_1(z,s)L_{1}(U^{n-1}(s),U^{n-1}(s-\tau))\big|^2ds\\
 &\leqslant&5\sup_{z\in[0,t]}\|\Phi(t)\|^2_\FF+5\|g_2\|^2_\infty t\int_{0}^{t}\EE \|K_{2}(U^{n-1}(s),U^{n-1}(s-\tau))\|^2_\FF ds\\
&&+\ 5\|g_1\|^2_\infty  t\int_0^t\EE \|K_{1}(U^{n-1}(s),U^{n-1}(s-\tau))\|^2_\FF ds\\
&&+\ 5\|h_2\|^2_\infty\int_0^t\EE \big|L_{2}(U^{n-1}(s),U^{n-1}(s-\tau))\big|^2ds\\
&&+\ 5\|h_1\|^2_\infty \int_0^t \EE\big|L_{1}(U^{n-1}(s),U^{n-1}(s-\tau))\big|^2ds.
  \end{eqnarray*}
Now, using hypothesis (H2) we get
\begin{eqnarray*}
&&\sup_{z\in[0,t]}\EE \|U^{n}(t)\|^2_\FF\\
 &\leqslant&5\sup_{z\in[0,t]}\|\Phi(t)\|^2_\FF+5C\big[t\|g_2\|^2_\infty +t\|g_1\|^2_\infty  +\|h_2\|^2_\infty+\|h_1\|^2_\infty\big]\\
&&\times\ \int_{0}^{t}\EE \big[1+\|U^{n-1}(s)\|^2_\FF +\|U^{n-1}(s-\tau))\|^2_\FF\big]ds\\
&\leqslant&M_1+M_2\int_{0}^{t}\EE \big[\|U^{n-1}(s)\|^2_\FF+\|U^{n-1}(s-\tau))\|^2_\FF\big]ds\\
&\leqslant&M_1+M_2\int_{0}^{t}\bigg[\sup_{z\in[0,s]}\EE \|U^{n-1}(z)\|^2_\FF+\sup_{z\in[0,s]}\EE \|U^{n-1}(z-\tau))\|^2_\FF\bigg]ds,
\end{eqnarray*}
where $M_1=5\|\Phi\|^{2}_{\II}+5C\tilde{T}[\tilde{T}\|g_2\|^2_\infty + \tilde{T}\|g_1\|^2_\infty+ \|h_2\|^2_\infty+ \|h_1\|^2_\infty]$ and $M_2=M_1/\tilde{T}$.

\noindent
Noticing that
$$
\int_0^t\sup_{z\in[0,s]}\EE \|U^{n-1}(z-\tau))\|^2_\FF ds\leqslant \tilde{T}\sup_{z\in[-\tau,0]}\|\Phi(z)\|^2_\FF+\int_0^t\sup_{z\in[0,s]}\EE \|U^{n-1}(z))\|^2_\FF ds,
$$
one can write
$$
\sup_{z\in[0,t]}\EE \|U^{n}(t)\|^2_\FF\leqslant M_3+M_4\int_{0}^{t}\sup_{z\in[0,s]}\EE \|U^{n-1}(z)\|^2_\FF ds,
$$
where $M_3,M_4$ are some new positive constants.

Let us observe that
for $m\geqslant 1$ we have
$$
\max_{1\leqslant n\leqslant m}\sup\limits_{z\in[0,t]}\EE \|U^{n}(t)\|^2_\FF\leqslant M_3+M_4\int_{0}^{t}\max_{1\leqslant n\leqslant m}\sup\limits_{z\in[0,s]}\EE \|U^{n-1}(t)\|^2_\FF ds.
$$
Further, inequality $\max\limits_{1\leqslant n\leqslant m}\sup\limits_{z\in[0,s]}\EE \|U^{n-1}(t)\|^2_\FF\leqslant \sup\limits_{z\in[-\tau,\tilde{T}]}\|\Phi(z)\|^2_\FF+\sup\limits_{z\in[0,s]}\EE \|U^{n}(t)\|^2_\FF$ leads us to
$$
\max\limits_{1\leqslant n\leqslant m}\sup\limits_{z\in[0,t]}\EE \|U^{n}(t)\|^2_\FF\leqslant M_5+M_4\int_{0}^{t}\max\limits_{1\leqslant n\leqslant m}\sup\limits_{z\in[0,s]}\EE \|U^{n}(t)\|^2_\FF ds,
$$
where $M_5$ is a certain new positive constant. 
Applying the Gronwall inequality to the above written inequality, we can infer that
$$
\max_{1\leqslant n\leqslant m}\sup\limits_{z\in[0,t]}\EE \|U^{n}(t)\|^2_\FF\leqslant M_5e^{M_4t}
$$
for every $t\in[0,\tilde{T}]$ and every $k\geqslant 1$.
Thus
$$
\max_{1\leqslant n\leqslant m}\sup\limits_{z\in[0,\tilde{T}]}\EE \|U^{n}(t)\|^2_\FF\leqslant M_5e^{M_4\tilde{T}}
$$
for every  $m\geqslant 1$, and this ends the proof.\qed

As a consequence of Lemma~\ref{OGR}, we obtain
$$
\EE\int_0^{\tilde{T}}\|U^n(t)\|^2_\FF dt
\leqslant \int_0^{\tilde{T}}\sup_{t\in[0,\tilde{T}]}\EE \|U^n(t)\|^2_\FF dt\leqslant M_5\tilde{T}e^{M_4\tilde{T}}<\infty
$$
 and this way the fuzzy stochastic process $U^n\Bigl|_{[0,\tilde{T}]\times\Omega}$ belongs to $\mathcal{L}^2([0,\tilde{T}]\times\Omega,\mathcal{N};\FF)$ for every $n\in\NN$.

The sequence of fuzzy stochastic processes $\{U^n\}$ has special significance. Based on it, the proof of the following theorem about the existence of a solution to Equation~\eqref{glowne} is constructed.

\begin{theorem}\label{existence and uniqueness theorem}
Let hypotheses (H0)-(H3) be satisfied. Then Equation~\eqref{glowne} has a unique solution $U$ on the interval $[-\tau,\tilde{T}]$.
\end{theorem}
\proof
First, we find an upper bound on the expression $\EE \sup_{z\in[0,t]}d^2_\infty (U^{1}(z),U^0(z))$, where $t$ is a fixed number in the interval $[0,\tilde{T}]$. Let us write
\begin{eqnarray*}
&&\EE \sup_{z\in[0,t]}d^2_\infty (U^{1}(z),U^0(z))\\
&= &\EE \sup_{z\in[0,t]}d^2_\infty \bigg(\bigg[\Phi(z)\oplus\int_{0}^{z}g_2(z,s)K_{2}(\Phi(s),\Phi(s-\tau))ds\bigg]\\
&&\ominus \int_{0}^{z}g_1(z,s)K_{1}\left(\Phi(s),\Phi(s-\tau)\right)ds\\
&&\oplus\bigg\langle\int_{0}^{z}h_2(z,s)L_{2}(\Phi(s),\Phi(s-\tau))dB_2(s)\\
&&-\int_{0}^{z}h_1(z,s)L_{1}(\Phi(s),\Phi(s-\tau))dB_1(s)\bigg\rangle,\Phi(z)\bigg)\\
&\leqslant&\EE\sup_{z\in[0,t]}\bigg[4 \bigg\|\int_{0}^{z}g_2(z,s)K_{2}(\Phi(s),\Phi(s-\tau))ds\bigg\|^2_\FF\\
&&+\ 4 \bigg\|\int_{0}^{z}g_1(z,s)K_{1}(\Phi(s),\Phi(s-\tau))ds\bigg\|^2_\FF\\
&&+\ 4 \bigg|\int_{0}^{z}h_2(z,s)L_{2}(\Phi(s),\Phi(s-\tau))dB_2(s)\bigg|^2\\
&&+\ 4 \bigg|\int_{0}^{z}h_1(z,s)L_{1}(\Phi(s),\Phi(s-\tau))dB_1(s)\bigg|^2\bigg].
\end{eqnarray*}
Applying Proposition 2.1 (iv), the Doob inequality and hypothesis (H2) we get
\begin{eqnarray*}
&&\sup_{z\in[0,t]}\EE d^2_\infty (U^{1}(z),U^0(z))\\
&\leqslant&4\|g_2\|^2_\infty t \EE\int_0^t \|K_{2}(\Phi(s),\Phi(s-\tau))\|^2_\FF ds\\
&&+\ 4\|g_1\|^2_\infty t\EE \int_0^t \|K_{1}(\Phi(s),\Phi(s-\tau))\|^2_\FF ds\\
&&+\ 16\|h_2\|^2_\infty \EE\int_0^t\big|L_{2}(\Phi(s),\Phi(s-\tau))\big|^2ds\\
&&+\ 16\|h_1\|^2_\infty \EE\int_0^t \big|L_{1}(\Phi(s),\Phi(s-\tau))\big|^2ds\\
&\leqslant&4C\big[\|g_2\|^2_\infty \tilde{T}+\|g_1\|^2_\infty  \tilde{T}+4\|h_2\|^2_\infty+4\|h_1\|^2_\infty\big]\\
&&\times\ \int_{0}^{t} \big[1+\|\Phi(s)\|^2_\FF+\|\Phi(s-\tau))\|^2_\FF\big]ds.
\end{eqnarray*}
This lead us to an estimation
$$
\sup_{z\in[0,t]}\EE d^2_\infty (U^{1}(z),U^0(z))\leqslant 4C\big[\|g_2\|^2_\infty \tilde{T}+\|g_1\|^2_\infty  \tilde{T}+4\|h_2\|^2_\infty+4\|h_1\|^2_\infty\big](1+2\sup\limits_{t\in[-\tau,\tilde{T}]}\|\Phi(t)\|^2_\FF)t.
$$
In the last expression, there is a constant number on the left side of the variable $t$. In the rest of the proof it will be symbolized by the letter $\xi$.

Now we will examine the upper bound of the expression $\EE \sup_{z\in[0,t]}d^2_\infty (U^{n+1}(z),U^{n}(z))$,  where $n\geqslant 1$ and $t\in\tilde{T}$. Notice that
\begin{eqnarray*}
&&\EE \sup_{z\in[0,t]}d^2_\infty (U^{n+1}(z),U^{n}(z))\\
 &\leqslant&4\|g_2\|^2_\infty t\EE\int_0^td^2_\infty (K_2(U^n(s),U^n(s-\tau)),K_2(U^{n-1}(s),U^{n-1}(s-\tau)))ds\\
&&+\ 4\|g_1\|^2_\infty  t\EE\int_0^td^2_\infty (K_1(U^n(s),U^n(s-\tau)),K_1(U^{n-1}(s),U^{n-1}(s-\tau)))ds\\
&&+\ 16\|h_2\|^2_\infty\EE\int_0^td^2_\infty (L_2(U^n(s),U^n(s-\tau)),L_2(U^{n-1}(s),U^{n-1}(s-\tau)))ds\\
&&+\ 16\|h_1\|^2_\infty\EE\int_0^td^2_\infty (L_1(U^n(s),U^n(s-\tau)),L_1(U^{n-1}(s),U^{n-1}(s-\tau)))ds.
\end{eqnarray*}
Invoking hypothesis (H3) we can write
\begin{eqnarray*}
&&\EE \sup_{z\in[0,t]}d^2_\infty (U^{n+1}(z),U^{n}(z))\\
&\leqslant&\zeta\EE\int_0^t\big[d^2_\infty (U^n(s),U^{n-1}(s))+d^2_\infty (U^n(s-\tau),U^{n-1}(s-\tau))\big]ds\\
&\leqslant&2\zeta\EE\int_0^td^2_\infty (U^n(s),U^{n-1}(s))ds\\
&\leqslant&2\zeta\int_0^t\EE\sup_{z\in[0,s]}d^2_\infty (U^n(z),U^{n-1}(z))ds,
\end{eqnarray*}
where $\zeta=4C\big[\|g_2\|^2_\infty \tilde{T}+ \|g_1\|^2_\infty \tilde{T}+ 4\|h_2\|^2_\infty+ 4\|h_2\|^2_\infty\big]$.
From the estimates we received, it can be concluded that
\begin{equation}\label{CZ}
\EE \sup_{z\in[0,t]}d^2_\infty (U^{n}(z),U^{n-1}(z))\leqslant \frac{\xi}{2\zeta}\cdot\frac{(2\zeta t)^{n}}{n!}
\end{equation}
for every $t\in[0,\tilde{T}]$ and every  $n\in\NN$.
Applying \eqref{CZ} to the Chebyshev inequality, we obtain
$$
P\bigg(\sup_{t\in[0,\tilde{T}]}d^2_\infty (U^{n}(s),U^{n-1}(s)\big)>\frac{1}{4^{n}}\bigg)\leqslant \frac{\xi}{2\zeta}\cdot\frac{(8\zeta\tilde{T})^n}{n!},
$$
Now by the Borel--Cantelli lemma, we arrive at
$$
P\bigg(\sup_{t\in[0,\tilde{T}]}d^2_\infty (U^{n}(s),U^{n-1}(s)\big)>\frac{1}{4^{n}} \text{ infinitely often}\bigg)=0
$$
and equivalently
$$
P\bigg(\sup_{t\in[0,\tilde{T}]}d_\infty(U^{n}(s),U^{n-1}(s)\big)>\frac{1}{2^{n}} \text{ infinitely often}\bigg)=0.
$$
Therefore, the sequence $\{U^{n}(\cdot,\omega)\}$ possesses a limit  fuzzy set-valued function $\widetilde{U}(\cdot,\omega):[0,\tilde{T}]\rightarrow\FF$ to which it converges on the interval $[0,\tilde{T}]$. The convergence happens for every $\omega\in\tilde{\Omega}$ such that $P(\tilde{\Omega})=1$. Now, we can define $U(t,\omega):=\Phi(t)$ for $(t,\omega)\in[-\tau,0]\times\Omega$, $U(t,\omega):=\widetilde{U}(t,\omega)$ for $(t,\omega)\in[0,\tilde{T}]\times\tilde{\Omega}$ and $U(t,\omega)$ as freely choosen fuzzy set in $\FF$ for $(t,\omega)\in[0,\tilde{T}]\times(\Omega\setminus\tilde{\Omega})$.  Since the following convergence 
$$
\lim_{n\rightarrow\infty}d_\infty(U^n(t,\omega),U(t,\omega))=0\text{ for }(t,\omega)\in[0,\tilde{T}]\times\tilde{\Omega},
$$
takes place, 
the mapping $U(t,\cdot):\Omega\to\FF$ is $\mathcal{A}_t$-measurable fuzzy random variable. Thus, $U\big|_{[0,\tilde{T}]\times\Omega}$ is a $d_\infty$-continuous $\{\mathcal{A}_t\}$-adapted fuzzy stochastic process, and this way it is a nonanticipating fuzzy stochastic process. Let us also observe that
\begin{eqnarray*}
\EE\int_0^{\tilde{T}}\|U(t)\|^2_\FF dt
&\leqslant& 
2\EE\int_0^{\tilde{T}}d^2_\infty (U(t),U^n(t))dt+2\EE\int_0^{\tilde{T}}\|U^n(t)\|^2_\FF dt\\
&\leqslant& 
2\tilde{T}\EE\sup_{t\in[0,\tilde{T}]}d^2_\infty (U(t),U^n(t))+2{\tilde{T}}\sup_{t\in[0,\tilde{T}]}\EE \|U^n(t)\|^2_\FF.
\end{eqnarray*}
In the last expression, the first term converges to zero and the second term is uniformly bounded in $n$ by Lemma~\ref{OGR}. Therefore  $\EE\int_0^{\tilde{T}}\|U(t)\|^2_\FF dt<\infty$ and $U\in\mathcal{L}^2([0,\tilde{T}]\times\Omega, \mathcal{N};\FF)$.
To show that $U$ is a solution to \eqref{glowne}, it is left to show that condition (4) from Definition~\ref{ROZ} is satisfied. 
To do this, we will first show that the expression \begin{eqnarray*}
&&
S_n:=\EE \sup_{z\in[0,\tilde{T}]}d^2_\infty \bigg(U^n(z),\bigg[\Phi(z)\oplus\int_{0}^{z}g_2(z,s)K_{2}(U(s),U(s-\tau))ds\bigg]\\
&&\ominus \int_{0}^{z}g_1(z,s)K_{1}\left(U(s),U(s-\tau)\right)ds\\
&&\oplus\bigg\langle\int_{0}^{z}h_2(z,s)L_{2}(U(s),U(s-\tau))dB_2(s)-\int_{0}^{z}h_1(z,s)L_{1}(U(s),U(s-\tau))dB_1(s)\bigg\rangle\bigg)
\end{eqnarray*} converges to zero as $n\rightarrow\infty$. Indeed,
\begin{eqnarray*}
S_n&\leqslant&4\EE \sup_{z\in[0,\tilde{T}]}d^2_\infty \bigg(\int_{0}^{z}g_2(z,s)K_{2}(U^{n-1}(s),U^{n-1}(s-\tau))ds,\\
&&\int_{0}^{z}g_2(z,s)K_{2}(U(s),U(s-\tau))ds\bigg)\\
&&+\ 4\EE \sup_{z\in[0,\tilde{T}]}d^2_\infty \bigg(\int_{0}^{z}g_1(z,s)K_{1}(U^{n-1}(s),U^{n-1}(s-\tau))ds,\\
&&\int_{0}^{z}g_1(z,s)K_{1}(U(s),U(s-\tau))ds\bigg)\\
&&+\ 4\EE \sup_{z\in[0,\tilde{T}]}\bigg|\int_{0}^{z}h_2(z,s)L_{2}(U^{n-1}(s),U^{n-1}(s-\tau))dB_2(s)\\
&&-\int_{0}^{z}h_2(z,s)L_{2}(U(s),U(s-\tau))dB_2(s)\bigg|^2\\
&&+4\EE \sup_{z\in[0,\tilde{T}]}\bigg|\int_{0}^{z}h_1(z,s)L_{1}(U^{n-1}(s),U^{n-1}(s-\tau))dB_1(s)\\
&&-\int_{0}^{z}h_1(z,s)L_{1}(U(s),U(s-\tau))dB_1(s)\bigg|^2\\
&\leqslant& \kappa\EE\int_0^{\tilde{T}}\big[d^2_\infty (U^{n-1}(s),U(s))+d^2_\infty (U^{n-1}(s-\tau),U(s-\tau))\big]ds\\
&\leqslant& \kappa\EE\int_0^{\tilde{T}}\big[\sup_{s\in[0,\tilde{T}]}d^2_\infty (U^{n-1}(s),U(s))+\sup_{s\in[0,\tilde{T}]}d^2_\infty (U^{n-1}(s-\tau),U(s-\tau))\big]ds\\
&\leqslant&\delta\EE\sup_{s\in[0,\tilde{T}]}d^2_\infty (U^{n-1}(s),U(s))\stackrel{n\rightarrow\infty}{\longrightarrow}0,
\end{eqnarray*}
where $\kappa=4C\big[\|g_2\|^2_\infty\tilde{T}+\|g_1\|^2\tilde{T}+4\|h_2\|^2_\infty+4\|h_1\|^2_\infty\big]$ and $\delta=2\kappa\tilde{T}$. 
Let us notice that
\begin{eqnarray*}
&&\EE \sup_{z\in[0,\tilde{T}]}d^2_\infty \bigg(U(z),\bigg[\Phi(z)\oplus\int_{0}^{z}g_2(z,s)K_{2}(U(s),U(s-\tau))ds\bigg]\\
&&\ominus \int_{0}^{z}g_1(z,s)K_{1}\left(U(s),U(s-\tau)\right)ds\\
&&\oplus\bigg\langle\int_{0}^{z}h_2(z,s)L_{2}(U(s),U(s-\tau))dB_2(s)-\int_{0}^{z}h_1(z,s)L_{1}(U(s),U(s-\tau))dB_1(s)\bigg\rangle\bigg)\\
&\leqslant&2\EE \sup_{z\in[0,\tilde{T}]}d^2_\infty (U(z),U^n(z))+2S_n
\end{eqnarray*}
and the right hand side of this inequality converges to zero which implies that the left hand side of the inequality equals zero. Hence we infer that
\begin{eqnarray*}
&&P\bigg( \sup_{z\in[0,\tilde{T}]}H\bigg(U(z),\bigg[\Phi(z)\oplus\int_{0}^{z}g_2(z,s)K_{2}(U(s),U(s-\tau))ds\bigg]\\
&&\ominus \int_{0}^{z}g_1(z,s)K_{1}\left(U(s),U(s-\tau)\right)ds\\
&&\oplus\bigg\langle\int_{0}^{z}h_2(z,s)L_{2}(U(s),U(s-\tau))dB_2(s)\\
&&-\int_{0}^{z}h_1(z,s)L_{1}(U(s),U(s-\tau))dB_1(s)\bigg\rangle=0\bigg)=1.
\end{eqnarray*}
This completes the justification that $U$ is a solution to \eqref{glowne}. Now we will show that $U$ is a unique solution. Let us say that a fuzzy stochastic process $V\colon[-\tau,\tilde{T}]\times\Omega\to\FF$ satisfies Definition~\ref{ROZ} too. Then for every $t\in[0,\tilde{T}]$ we obtain
\begin{eqnarray*}
&&\EE \sup_{z\in[0,t]}d^2_\infty (U(z),V(z))\\
&\leqslant& \kappa\EE\int_0^t\big[\sup_{z\in[0,s]}d^2_\infty (U(z),V(z))+\sup_{z\in[0,s]}d^2_\infty (U(z-\tau),V(z-\tau))\big]ds\\
&\leqslant&2\kappa\int_0^t\EE\sup_{z\in[0,s]}d^2_\infty (U(z),V(z))ds.
\end{eqnarray*}
With the help of Gronwall's inequality we get
$$
\EE \sup_{z\in[0,t]}d^2_\infty (U(z),V(z))\leqslant 0 \text{ for }t\in[0,\tilde{T}]
$$
and this immediately implies that $\EE \sup\limits_{z\in[0,t]}d^2_\infty (U(z),V(z))=0 \text{ for }t\in[0,\tilde{T}]
$. Hence
$$
P\bigg(\sup\limits_{z\in[0,\tilde{T}]}d_\infty (U(z),V(z))=0\bigg)=1
$$
and according to the Definition~\ref{JED} this means that the solution is the only one.
\qed

Now we consider the equation
\begin{eqnarray}\label{INI}
V(t)&\oplus& \int_0^tg_1(t,s)K_1(V(s),V(s-\tau))ds\oplus \Bigl\langle\int_0^th_1(t,s)L_1(V(s),V(s-\tau))dB_1(s)\Bigr\rangle\nonumber\\
&=&\Psi(t)\oplus\int_0^tg_2(t,s)K_2(V(s),V(s-\tau))ds\nonumber\\
&&\oplus\ \Bigl\langle\int_0^th_2(t,s)L_2(V(s),V(s-\tau))dB_2(s)\Bigr\rangle,\quad t\in[0,T],\\
V(t)&=&\Psi(t),\quad t\in[-\tau,0]\nonumber
\end{eqnarray}
which differs from \eqref{glowne} in its initial value. Now this initial value is a function $\Psi$. The intention of further research on symmetric fuzzy
stochastic Volterra integral equations with retardation is to justify the fact that if $\Psi$ is only slightly different from the original initial value $\Phi$, it does not have a major impact on the solution. In fact, the solutions are not far apart.

\begin{theorem} Assume that (H0) is satisfied by $\Phi$ and $\Psi$. Suppose that condition (H1) is met. Let $U:[-\tau,\tilde{T}]\times\Omega\to\FF$ and  $V:[-\tau,\tilde{T}]\times\Omega\to\FF$ be the solutions to \eqref{glowne} and \eqref{INI}, respectively.
Then  $$\EE\sup_{t\in[-\tau,\tilde{T}]}d^2_\infty (U(t),V(t))\leqslant M\sup_{[-\tau,\tilde{T}]}d^2_\infty (\Phi(t),\Psi(t)),$$
where $M$ is a positive number.
\end{theorem}
\proof
First, let us note that for $t\in[-\tau,0]$ it holds $V(t)=\Psi(t)$ and $U(t)=\Phi(t)$. Therefore $$\EE\sup\limits_{t\in[-\tau,0]}d^2_\infty (U(t),V(t))=\sup\limits_{t\in[-\tau,0]}d^2_\infty (\Phi(t),\Psi(t)).$$

\noindent Now let us consider for $t\in[0,\tilde{T}]$. In this case
\begin{eqnarray*}
&&\hspace{-10mm}\EE \sup_{z\in[0,t]}d^2_\infty (U(z),V(z))\\
&=&\EE\sup_{z\in[0,t]}d^2_\infty \Bigl(\Bigl[\Phi(t)\oplus\int_0^tg_2(t,s)K_2(U(s),U(s-\tau))ds\Bigr]\\
&&\ominus \int_0^tg_1(t,s)K_1(U(s),U(s-\tau))ds\\
&&\oplus    \Bigl\langle\int_0^th_2(t,s)L_2(U(s),U(s-\tau))dB_2(s)-\int_0^th_1(t,s)L_1(U(s),U(s-\tau))dB_1(s)\Bigr\rangle,\\
&& \Bigl[\Psi(t)\oplus\int_0^tg_2(t,s)K_2(V(s),V(s-\tau))ds\Bigr]\ominus \int_0^tg_1(t,s)K_1(V(s),V(s-\tau))ds\\
&&\oplus   \Bigl\langle\int_0^th_2(t,s)L_2(V(s),V(s-\tau))dB_2(s)-\int_0^th_1(t,s)L_1(V(s),V(s-\tau))dB_1(s)\Bigr\rangle\Bigr)\\
&\leqslant&5\sup_{z\in[0,t]}d^2_\infty (\Phi(z),\Psi(z))\\
&&+\ 5\EE\sup_{z\in[0,t]}d^2_\infty (\int_0^zg_2(z,s)K_2(U(s),U(s-\tau))ds,\int_0^zg_2(z,s)K_2(V(s),V(s-\tau))ds)\\
&&+\ 5\EE\sup_{z\in[0,t]}d^2_\infty (\int_0^zg_1(z,s)K_1(U(s),U(s-\tau))ds,\int_0^zg_1(z,s)K_1(V(s),V(s-\tau))ds)\\
&&+\ 5\EE\sup_{z\in[0,t]}\Bigl|\int_0^z[h_2(z,s)L_2(U(s),U(s-\tau))-h_2(z,s)L_2(V(s),V(s-\tau))]dB_2(s)\Bigr|^2\\
&&+\ 5\EE\sup_{z\in[0,t]}\Bigl|\int_0^z[h_1(z,s)L_1(U(s),U(s-\tau))-h_1(z,s)L_1(V(s),V(s-\tau))]dB_1(s)\Bigr|^2.
\end{eqnarray*}
Through application of Proposition~\ref{PRO} (iv) and the Doob maximal inequality we obtain 
\begin{eqnarray*}
&&\hspace{-10mm}\EE \sup_{z\in[0,t]}d^2_\infty (U(z),V(z))\\
&\leqslant&5\sup_{z\in[0,t]}d^2_\infty (\Phi(z),\Psi(z))\\
&&+\ 5\|g_2\|^2t\EE\int_0^td^2_\infty (K_2(U(s),U(s-\tau)),K_2(V(s),V(s-\tau)))ds\\
&&+\ 5\|g_1\|^2t\EE\int_0^td^2_\infty (K_1(U(s),U(s-\tau)),K_1(V(s),V(s-\tau)))ds\\
&&+\ 20\|h_2\|^2\EE\int_0^t\Bigl|L_2(U(s),U(s-\tau))-L_2(V(s),V(s-\tau))\Bigr|^2ds\\
&&+\ 20\|h_1\|^2\EE\int_0^t\Bigl|L_1(U(s),U(s-\tau))-L_1(V(s),V(s-\tau))\Bigr|^2ds.
\end{eqnarray*}
Now by referring to (H1) we are led to 
\begin{eqnarray*}
&&\hspace{-10mm}\EE \sup_{z\in[0,t]}d^2_\infty (U(z),V(z))\\
&\leqslant&5\sup_{z\in[0,t]}d^2_\infty (\Phi(z),\Psi(z))+\gamma \EE \int_0^t\Bigl(d^2_\infty (U(s),V(s))+d^2_\infty (U(s-\tau),V(s-\tau))\Bigr)ds,
\end{eqnarray*}
where $\gamma=5\|g_2\|_\infty^2\tilde{T}C+ 5\|g_1\|_\infty^2\tilde{T}C+ 20\|h_2\|_\infty^2C+20\|h_1\|_\infty^2C$. A further step is to notice that
\begin{eqnarray*}
\EE \sup_{z\in[0,t]}d^2_\infty (U(z),V(z))&\leqslant &5\sup_{z\in[0,t]}d^2_\infty (\Phi(z),\Psi(z))+ 2\gamma\EE \int_0^td^2_\infty (U(s),V(s))ds\\
&\leqslant &5\sup_{z\in[0,t]}d^2_\infty (\Phi(z),\Psi(z))+ 2\gamma \int_0^t\EE\sup_{z\in[0,s]}d^2_\infty (U(z),V(z))ds.
\end{eqnarray*}
Referring to Gronwall's inequality, we can write
\begin{eqnarray*}
\EE\sup_{z\in[0,t]}d^2_\infty (U(z),V(z))&\leqslant&5\sup_{z\in[0,t]}d^2_\infty (\Phi(z),\Psi(z))\exp\{2\gamma t\}\quad\text{for every}\quad t\in[0,\tilde{T}].
\end{eqnarray*}
Therefore
$$
\EE\sup_{z\in[0,\tilde{T}]}d^2_\infty (U(z),V(z))\leqslant 5\exp\{2\gamma 
\tilde{T}\}\sup_{z\in[0,\tilde{T}]}d^2_\infty (\Phi(z),\Psi(z))
$$
and this leads to the end of the proof.
\qed

From the above theorem we obtain the property of continuous dependence of the solution $U$ on the initial value. Indeed, if we consider a sequence of equations of type \eqref{INI} with initial values $\Psi_n$ and solutions $V_n$ and property $\sup_{t\in[-\tau,T]}d_\infty(\Phi(t),\Psi_n(t))\rightarrow 0$, as $n\rightarrow\infty$, would be satisfied, then due to the above theorem we will obtain good behavior of the solutions, i.e.\
$$
\EE\sup_{z\in[0,\tilde{T}]}d^2_\infty (U(z),V_n(z))\rightarrow 0\text{ as }n\rightarrow\infty.
$$

In the next part of this section, we will argue that we will obtain a similar property when we consider a sequence of equations with other nonlinearities and other kernels, but those that in a sense converge to the original data in \eqref{glowne}. Therefore we consider \eqref{glowne} where the kernels are $g_1,g_2,h_1,h_2$ and nonlinearities $K_1,K_2,L_1,L_2$ and also consider other kernels $g_1^n,g_2^n,h_1^n,h_2^n$ and other nonlinearities $K_1^n,K_2^n,L_1^n,L_2^n$, where $n=1,2,\ldots$ in a sequence of equations
\begin{eqnarray}\label{NON}
V^n(t)&\oplus&\int_0^tg^n_1(t,s)K_1^n(V^n(s),V^n(s-\tau))ds\oplus\Bigl\langle
\int_0^th_1^n(t,s)L_1^n(V^n(s),V^n(s-\tau))dB_1(s)\Bigr\rangle\nonumber\\
&=&\Phi(t)\oplus\int_0^tg^n_2(t,s)K_2^n(V^n(s),V^n(s-\tau))ds\nonumber\\
&&\oplus\Bigl\langle\int_0^th_2^n(t,s)L_2^n(V^n(s),V^n(s-\tau))dB_2(s)\Bigr\rangle,\quad t\in[0,T],\\
V^n(t)&=&\Phi(t),\quad t\in[-\tau,0]\text{ for }n=1,2,\ldots\nonumber
\end{eqnarray}
 \begin{theorem} Let $U:[-\tau,\tilde{T}]\times\Omega\to\FF$ and $V^n:[-\tau,\tilde{T}]\times\Omega\to\FF$ be the solutions to \eqref{glowne} and \eqref{NON}, respectively. Assume that   $\max\limits_n\|g_1^n\|_\infty+\max\limits_n\|h_1^n\|_\infty+\max\limits_n\|g_2^n\|_\infty+\max\limits_n\|h_2^n\|_\infty<\infty$ and condition (H1) is satisfied.
If for every $t\in[0,T]$ and for every $U,V\in\FF$ it holds
\begin{eqnarray}\label{zaq}
&&\int_0^Td^2_\infty (g_1^n(t,s)K_1^n(U,V),g_1(t,s)K_1(U,V))ds\nonumber\\
&&
+\ \int_0^Td^2_\infty (g_2^n(t,s)K_2^n(U,V),g_2(t,s)K_2(U,V))ds\nonumber\\
&&
+\ \int_0^T\Bigl| h_1^n(t,s)L_1^n(U,V)- h_1(t,s)L_1(U,V)\Bigr|^2ds\nonumber\\
&&
+\ \int_0^T\Bigl| h_2^n(t,s)L_2^n(U,V)- h_2(t,s)L_2(U,V)\Bigr|^2ds\longrightarrow 0 \text{ as } n\rightarrow\infty,
\end{eqnarray}
then $$\EE d^2_\infty (U(t),V^n(t))\longrightarrow 0 \text{ as }n\rightarrow\infty \text{ for every }t\in[-\tau,T].$$
\end{theorem}
\proof
Since $V^n(t)=\Phi(t)$ and $U(t)=\Phi(t)$  for $t\in[-\tau,0]$,  we immediately have that 
$\EE d^2_\infty (U(t),V^n(t))=0$ for $t\in[-\tau,0]$.

For $t\in[0,T]$ let us perceive
\begin{eqnarray*}
&&\hspace{-10mm}\EE d^2_\infty (U(t),V^n(t))\\
&\leqslant&4\EE d^2_\infty (\int_0^tg_2(t,s)K_2(U(s),U(s-\tau))ds,\int_0^tg_2^n(t,s)K_2^n(V^n(s),V^n(s-\tau))ds)\\
&&+\ 4\EE d^2_\infty (\int_0^tg_1(t,s)K_1(U(s),U(s-\tau))ds,\int_0^tg_1^n(t,s)K_1^n(V^n(s),V^n(s-\tau))ds)\\
&&+\ 4\EE \Bigl|\int_0^th_2(t,s)L_2(U(s),U(s-\tau))dB_2(s)\\
&&-\int_0^th_2^n(t,s)L_2^n(V^n(s),V^n(s-\tau))dB_2(s)\Bigr|^2\\
&&+\ 4\EE \Bigl|\int_0^th_1(t,s)L_1(U(s),U(s-\tau))dB_1(s)\\
&&-\int_0^th_1^n(t,s)L_1^n(V^n(s),V^n(s-\tau))dB_1(s)\Bigr|^2.
\end{eqnarray*}
Continuing further
\begin{eqnarray*}
&&\hspace{-10mm}\EE d^2_\infty (U(t),V^n(t))\\
&\leqslant&8\EE d^2_\infty (\int_0^tg_2(t,s)K_2(U(s),U(s-\tau))ds,\int_0^tg_2^n(t,s)K_2^n(U(s),U(s-\tau))ds)\\
&&+\ 8\EE d^2_\infty (\int_0^tg_2^n(t,s)K_2^n(U(s),U(s-\tau))ds,\int_0^tg_2^n(t,s)K_2^n(V^n(s),V^n(s-\tau))ds)\\
&&+\ 8\EE d^2_\infty (\int_0^tg_1(t,s)K_1(U(s),U(s-\tau))ds,\int_0^tg_1^n(t,s)K_1^n(U(s),U(s-\tau))ds)\\
&&+\ 8\EE d^2_\infty (\int_0^tg_1^n(t,s)K_1^n(U(s),U(s-\tau))ds,\int_0^tg_1^n(t,s)K_1^n(V^n(s),V^n(s-\tau))ds)\\
&&+\ 8\EE \Bigl|\int_0^th_2(t,s)L_2(U(s),U(s-\tau))dB_2(s)\\
&&-\ \int_0^th_2^n(t,s)L_2^n(U(s),U(s-\tau))dB_2(s)\Bigr|^2\\
&&+\ 8\EE \Bigl|\int_0^th_2^n(t,s)L_2^n(U(s),U(s-\tau))dB_2(s)\\
&&-\ \int_0^th_2^n(t,s)L_2^n(V^n(s),V^n(s-\tau))dB_2(s)\Bigr|^2\\
&&+\ 8\EE \Bigl|\int_0^th_1(t,s)L_1(U(s),U(s-\tau))dB_1(s)\\
&&-\ \int_0^th_1^n(t,s)L_1^n(U(s),U(s-\tau))dB_1(s)\Bigr|^2\\
&&+\ 8\EE \Bigl|\int_0^th_1^n(t,s)L_1^n(U(s),U(s-\tau))dB_1(s)\\
&&-\ \int_0^th_1^n(t,s)L_1^n(V^n(s),V^n(s-\tau))dB_1(s)\Bigr|^2.
\end{eqnarray*}
Due to Proposition~\ref{PRO} (iii) and the It\^o isometry 
\begin{eqnarray*}
&&\hspace{-10mm}\EE d^2_\infty (U(t),V^n(t))\\
&\leqslant&8\tilde{T}\EE \int_0^{\tilde{T}} d^2_\infty (g_2(t,s)K_2(U(s),U(s-\tau)),g_2^n(t,s)K_2^n(U(s),U(s-\tau)))ds\\
&&+\ 8\tilde{T}\EE \int_0^{\tilde{T}} d^2_\infty (g_1(t,s)K_1(U(s),U(s-\tau)),g_1^n(t,s)K_1^n(U(s),U(s-\tau)))ds\\
&&+\ 8\EE \int_0^{\tilde{T}}\Bigl|h_2(t,s)L_2(U(s),U(s-\tau))-h_2^n(t,s)L_2^n(U(s),U(s-\tau))\Bigr|^2ds\\
&&+\ 8\EE \int_0^{\tilde{T}}\Bigl|h_1(t,s)L_1(U(s),U(s-\tau))-h_1^n(t,s)L_1^n(U(s),U(s-\tau))\Bigr|^2ds\\
&&+\ 8\|g_2^n\|^2_\infty\tilde{T}\EE \int_0^t d^2_\infty (K_2^n(U(s),U(s-\tau)),K_2^n(V^n(s),V^n(s-\tau)))ds\\
&&+\ 8\|g_1^n\|^2_\infty\tilde{T}\EE \int_0^t d^2_\infty (K_1^n(U(s),U(s-\tau)),K_1^n(V^n(s),V^n(s-\tau)))ds\\
&&+\ 8\|h_2^n\|^2_\infty\EE \int_0^t\Bigl|L_2^n(U(s),U(s-\tau))-L_2^n(V^n(s),V^n(s-\tau))\Bigr|^2ds\\
&&+\ 8\|h_1^n\|^2_\infty\EE \int_0^t\Bigl|L_1^n(U(s),U(s-\tau))-L_1^n(V^n(s),V^n(s-\tau))\Bigr|^2ds.
\end{eqnarray*}
The group of the first four terms on the right side of the above inequality, let us denote it by $\delta_n$, converges to zero as $n\rightarrow \infty$ and this is due to the Lebesgue dominated convergence theorem and assumption (H1). Let us denote the group of the next four components on the right side of the above inequality by $\theta_n(t)$. Using (H1) and the Fubini theorem we arrive at
\begin{eqnarray*}
\theta_n(t)&\leqslant&\eta \int_0^t\EE \bigl(d^2_\infty (U(s),V^n(s))+d^2_\infty (U(s-\tau),V^n(s-\tau))\bigr)ds\\
&\leqslant&2\eta \int_0^t\EE d^2_\infty (U(s),V^n(s))ds,
\end{eqnarray*}
where $\eta= 8C\Bigl(\max_n\|g_2^n\|^2_\infty+ \max_n\|g_1^n\|^2_\infty\tilde{T}+ \max_n\|h_2^n\|^2_\infty+\max_n\|h_1^n\|^2_\infty\Bigr)$.
Hence
\begin{eqnarray*}
\EE d^2_\infty (U(t),V^n(t))&\leqslant&\delta_n+2\eta \int_0^t\EE d^2_\infty (U(s),V^n(s))ds\quad
\text{for every}\quad t\in[0,\tilde{T}].
\end{eqnarray*}
Applying the Gronwall inequality we obtain
\begin{eqnarray*}\label{edc}
\EE d^2_\infty (V_n(t),U(t))&\leqslant&\delta_n\exp\{2\eta t\}\quad
\text{for every}\quad t\in[0,\tilde{T}],
\end{eqnarray*}
Due to the fact that $\delta_n\stackrel{n\rightarrow\infty}{\longrightarrow} 0$ one gets
$$
\EE d^2_\infty (U(t),V_n(t))\stackrel{n\rightarrow\infty}{\longrightarrow} 0\quad \text{for every}\quad t\in[0,\tilde{T}].
$$
This ends the proof.
\qed

\section{Concluding remarks}
In this paper, symmetric fuzzy stochastic Volterra integral equations with constant retardation have been studied. The equations are called symmetric because in their formula Lebesgue and It\^o-type integrals appear symmetrically on both sides of the equation. This formulation of equations is appropriate for studies with fuzzy sets. The first result presented in the paper is the theorem on the existence and uniqueness of the solution, which was obtained with the requirements of the Lipschitz drift and diffusion coefficients. It is also reasoned that the solutions depend continuously on the equation parameters such as the initial value, nonlinearities and kernels. This confirms that the solution to the equation under consideration is stable in the sense that by changing the parameters of the equation a little, it is not possible for the solution to change much. Such property of solutions is naturally desirable.

\end{document}